\documentclass[11pt]{amsart}

\usepackage{amssymb,amsmath,latexsym, amscd}
\input cyracc.def

\newtheorem{theorem}{Theorem}[section]
\newtheorem{lemma}[theorem]{Lemma}
\newtheorem{proposition}[theorem]{Proposition}
\newtheorem{definition}[theorem]{Definition}

\newtheorem{corollary}[theorem]{Corollary}

\newtheorem{question}[theorem]{Question}

\def\s{\sigma}

\begin{document}

\title{Isolated points in the space of left orderings of a group}

\date{\today}

\author[Adam Clay]{Adam Clay}
\address{Department of Mathematics\\
University of British Columbia \\
Vancouver \\
BC Canada V6T 1Z2} \email{aclay@math.ubc.ca}
\urladdr{http://www.math.ubc.ca/~aclay/} \maketitle

\begin{abstract}
Let $G$ be a left orderable group and $LO(G)$ the space of all left orderings.   We investigate the circumstances under which a left ordering $<$ of $G$ can correspond to an isolated point in $LO(G)$, in particular we extend the main result of \cite{NF07} to the case of uncountable groups.   With minor technical restrictions on the group $G$, we find that no dense left ordering is isolated in $LO(G)$, and that the closure of the set of all dense left orderings of $G$ yields a dense $G_{\delta}$ set within a Cantor set of left orderings in $LO(G)$.  Lastly, we show that certain conditions on a discrete left ordering of $G$ can guarantee that it is not isolated in $LO(G)$, and we illustrate these ideas using the Dehornoy ordering of the braid groups.
\end{abstract}

\section{The space of left orderings of a group}

A group $G$ is said to be left-orderable if there exists a strict total ordering $<$ of its elements such that $g<h \Rightarrow fg < fh$ for all $f, g, h \in G$.  Given a left-orderable group $G$ with ordering $<$, we can identify the left ordering $<$ of $G$ with its positive cone $P =\{ g \in G | g>1 \}$, the set of all positive elements.  The positive cone $P$ of a left ordering of a group $G$ satisfies the following two properties:
\begin{enumerate}
\item If $g, h \in P$ then $gh \in P$.
\item For all $g \in G$, exactly one of $g \in P, g^{-1} \in P$, or $g =1$ holds.
\end{enumerate}
Conversely, given a semigroup $P \subset G$ satisfying the above two properties, we can order the elements of $G$ by specifying that $g<h$ if and only if $g^{-1}h \in P$.  

A left ordering $<$ of $G$ is said to be a Conradian ordering if whenever $g, h >1$, then there exists $n \in \mathbb{N}$ such that $g<hg^n$.  Lastly, a left ordering of a group $G$ is said to be a bi-ordering if the ordering is also invariant under multiplication from the right, namely $g<h \Rightarrow gf<hf$ for all $f, g, h \in G$.  It should be noted that the positive cone $P \subset G$ of a bi-ordering also satsifies the additional property:
\begin{enumerate}
\setcounter{enumi}{2}
\item For all $g \in G$, we have $gPg^{-1} =P$.
\end{enumerate}
Analogous to the case of left orderings, a semigroup $P \subset G$ satisfying properties (1)--(3) defines a bi-ordering of $G$.

  We can then consider the set $LO(G) \subset 2^G$ of all positive cones in $G$, a space first defined in \cite{AS04}.  As there is a one-to-one correspondence between left orderings of $G$ and positive cones in $G$, it is natural to describe $LO(G)$ as the space of all left orderings of $G$.  The space $LO(G)$ is endowed with the subspace topology arising from the product topology on $2^G$, with a subbasis for the topology on $LO(G)$ being formed by the open sets $U_g = \{ P \in LO(G) | g \in P \}$.  Note that $LO(G)$ comes equipped with a natural $G$-action: given an element $g \in G$, the positive cone $P$ is sent by $g$ to its conjugate $gPg^{-1}$.  Therefore, given a left ordering $<$ of $G$ with positive cone $P$, we can create new left orderings of $G$ by conjugating the corresponding positive cone $P$ by different elements of $G$.

One can check that $LO(G)$ is a closed subset in $2^G$, and from Tychonoff's Theorem we know that $2^G$ is compact, so that  $LO(G)$ itself must be a compact space.   With this setup, it is also easy to see that $LO(G)$ is a totally disconnected Hausdorff space, and in \cite{AS04} it is shown that whenever $G$ is countable, the topology on $LO(G)$ in fact arises from a very natural metric.  Thus we arrive at:

\begin{theorem}(Sikora, \cite{AS04}) 
Let $G$ be a countable group.  Then the space $LO(G)$ is a compact, totally disconnected Hausdorff metric space.  If $LO(G)$ also contains no isolated points, then $LO(G)$ is homeomorphic to the Cantor set. 
\end{theorem}

Given a group $G$, we would therefore like to address the existence of isolated points in the space $LO(G)$, as a first step towards understanding the structure of $LO(G)$.

Recall that a subgroup $C$ of a left-ordered group $G$ is called convex (with respect to the ordering $<$) if whenever $f, h \in C$ and $g \in G$, the implication $f<g<h \Rightarrow g \in C$ holds. For example, it is easy to check that the subgroup $C$ in Proposition \ref{prop:LO} is convex in the ordering contructed on $G$.

Following \cite{NF07}, we define the Conradian soul $C_<(G)$ in a left ordered group $G$ with ordering $<$ to be the largest convex subgroup $C \subset G$ such that the restriction of $<$ to $C$ is a Conradian ordering.  Similarly, we use the notation $B_<(G)$ to denote the largest convex subgroup $C \subset G$ such that the restriction of $<$ to $C$ is a bi-ordering.  Note that we always have $B_<(G) \subset C_<(G)$, since all bi-orderings are also Conradian orderings.

Using this notation, the main result of \cite{NF07}, which we will extend here to the case of uncountable groups, can be stated as follows.

\begin{theorem}
\label{thm:BCS}
Let $G$ be a group, and let $P \in LO(G)$ be an isolated point with associated ordering $<$ of $G$.  Then $B_<(G)$ is abelian of rank one, and $C_<(G)$ is non-trivial and admits only finitely many left orderings.
\end{theorem}

Note that Theorem \ref{thm:BCS} is proven for the case of countable groups in \cite{NF07}, although the dynamical approach used therein is entirely different than our approach, and does not generalize to the case of uncountable groups. 

Finally, recall that a left ordering of a group $G$ is dense if whenever $g<h$, then there exists $f \in G$ such that $g<f<h$.   If a left ordering $<$ of $G$ is not dense, then it is discrete, meaning that in the ordering $<$ of $G$ there is a least positive element $ \epsilon >1$.  We explore the structure of $LO(G)$ by considering the cases of dense and discrete left orderings separately, and we will find:

\begin{theorem}
\label{th:dense}
Let $Z \subset LO(G)$ denote the set of all dense left orderings of a countable group $G$, and suppose that all rank one abelian subgroups of $G$ are isomorphic
to $\mathbb{Z}$. Then if $Z$ is non-empty, its closure $\bar{Z}$ is homeomorphic to the Cantor set, and the set $Z$ is a $G_{\delta}$ set that is dense in $\bar{Z}$.
\end{theorem}

In the case of abelian groups, our result will be slightly stronger than Theorem \ref{th:dense}.  Specfically, in the case that $G$ is countable and abelian, we will show that $\bar{Z} = LO(G)$.

\textbf{Acknowledgments.} The author would like to thank Dale Rolfsen, Andr\'{e}s Navas and Crist\'{o}bal Rivas for many useful discussions and comments regarding earlier drafts of this paper.

\section{The case of Conradian orderings}

We first review known results concerning Conradian orderings, and consider also the case of bi-orderings.  Note that the results of this section concerning $C_<(G)$ appear in \cite{NF07}, and rely on the following difficult lemma (\cite{NF07} Lemma 4.4), the bulk of which appeared first in \cite{PL06}, and partially in \cite{KME96}.

\begin{lemma}
Suppose that $P$ is the positive cone of a Conradian ordering of the group $G$, and that there is exactly one proper, nontrivial convex subgroup $C \subset G$.  Further suppose that both $C$ and $G/C$ are rank one abelian groups.  If $P$ is isolated in $LO(G)$, then $G$ is not biorderable. \label{lem:2convex}
\end{lemma}

The next two theorems require the following work of Tararin (\cite{KME96}, Theorem 5.2.1).  Recall that a group $G$ admits a finite rational series if 
\[ 1 =G_0 \lhd G_1 \lhd \cdots \lhd G_n = G \]
is a finite normal series with all quotients $G_{i+1}/G_i$ rank one abelian.

\begin{theorem}
\label{thm:tararin}
Let $G$ be a left-ordered group.   
\begin{enumerate}
\item If $LO(G)$ is finite, then $G$ has a finite rational series.
\item Suppose that $G$ has a finite rational series.  Then $LO(G)$ is finite if and only if $G_i \lhd G$ for all $i$, and none of the quotients $G_{i+2}/G_i$ are bi-orderable.  Furthermore, in this case the rational series is unique, and for every left ordering of $G$, the convex subgroups are precisely $G_0, G_1, \cdots, G_n$.
\end{enumerate}
\end{theorem}

\begin{theorem} [\cite{NF07} Proposition 4.1]
\label{thm:consolated}
Suppose that $P$ is the positive cone of a Conradian ordering of $G$.  Then $P$ is not an isolated point in the space $LO(G)$, unless $LO(G)$ is finite.
\end{theorem}

\begin{theorem}
\label{thm:bisolated}
Suppose that $P$ is the positive cone of a bi-ordering of $G$.  Then $P$ is not isolated in $LO(G)$ unless $G$ is rank 1 abelian.
\end{theorem}
\begin{proof}
In the case that $G$ is bi-ordered by the ordering $<$ associated to $P$, we have $C_<(G) = G$.  From Theorem \ref{thm:consolated}, it follows that $G$ itself must have only finitely many left orderings if the bi-ordering $<$ is to have a positive cone that is isolated in $LO(G)$.  However, by the work of Tararin, we see that no group $G$ admitting only finitely many left orders is bi-orderable, except in the case that $G$ is rank one abelian.
\end{proof}

\section{Isolated points}

When trying to determine which points in $LO(G)$ are isolated, the conjugation action on $LO(G)$ is a useful tool in approximating a given positive cone. Aside from conjugation of a given ordering, there is a second natural way to make new left orderings of $G$, as follows.
\begin{proposition}
\label{prop:LO}
Suppose that $C$ is a left-orderable subgroup of $G$ with ordering $\prec$.  Suppose also that the left cosets of $C$ can be ordered in a way compatible with group multiplication from the left, namely $aC \prec' bC \Rightarrow caC \prec' cbC$ for all $a, b, c \in G$. Then a left ordering $<$ can be defined on $G$ by specifying a positive cone as follows:  An element $g \in G$ satisfies $1<g$ if $g \in C$ and $1 \prec g$, or if $g \notin C$ and $C \prec' gC$.
\end{proposition}
The proof is a simple check.  This proposition allows us to change any left ordering of a group $G$ on a specified convex subgroup $C$:  If $C \subset G$ is convex in the left ordering $<$, then convexity allows us to unambiguously define a left-invariant ordering $\prec$ of the cosets $\{gC | g \in G\}$.  We may then choose a left ordering of $C$ different from $<$, and extend it to a left ordering of $G$ by using the ordering $\prec$ of the cosets, and applying Proposition \ref{prop:LO}.

Next we observe some simple lemmas.

\begin{lemma}
Suppose $P \subset G$ and that $C$ is a convex subgroup of $G$.  Then if $P_C = P \cap C$ is not an isolated point in $LO(C)$, $P$ is not an isolated point in $LO(G)$.
\label{lem:congrp}
\end{lemma}
\begin{proof}
Suppose that 
\[ P \in \bigcap_{i=1}^m U_{g_i} ,\]
and suppose also that we have numbered the elements $g_i$ so that $g_i \in C$ for $i \leq k$ (possibly $k =0$, in the case that no $g_i$ lies in $C$).  Now in $LO(C)$, we have that 
\[ P_C \in \bigcap_{i=1}^k U_{g_i}, \]
and since $P_C$ is not an isolated point, we can choose $P_C' \in \bigcap_{i=1}^k U_{g_i}$, with $P_C' \neq P_C$.  

We can now construct a positive cone $P' \neq P$ on $G$ as follows:  Given $g \in G$, $g \in P'$ if $g \in C$ and $g \in P_C'$, or if $g \notin C$ and $g \in P$. 

The positive cone $P'$ is different from $P$, since $P$ and $P'$ disagree on $C$, and by construction, $P' \in \bigcap_{i=1}^m U_{g_i}$.  It follows that $P$ is not isolated.
\end{proof}

\begin{lemma}
Suppose $P \subset G$ and that $C$ is a normal, convex subgroup of $G$.  Let $P'$ denote the positive cone of the ordering inherited by the quotient $G/C$. 
If $P'$ is not an isolated point in $LO(G/C)$, $P$ is not an isolated point in $LO(G)$.
\label{lem:quotient}
\end{lemma}
The proof is routine.

\begin{lemma} Let $G$ be a left ordered group with ordering $<$, whose positive
cone we denote as $P$.
\label{lem:H} Then the subgroup
\[ stab(P) = \{ g\in G : gPg^{-1} = P \} \]
is bi-ordered by the restriction of $<$
to $H = stab(P)$.
\end{lemma}
\begin{proof}

To see that the restriction of $<$ is a bi-ordering, consider its
positive cone $P_H = P \cap H$.  If $g \in P_H$ and $h \in H$,
then \begin{itemize} \item $hgh^{-1} \in H$ since $H$ is a
subgroup, and \item $hgh^{-1} \in P$ since, by definition, every
element of $H$ fixes the positive cone $P$ under conjugation.
\end{itemize}
Therefore $H$ is bi-ordered.
\end{proof}

The main difficulty in characterizing the Conradian soul of an isolated point in $LO(G)$ is in showing that the Conradian soul is necessarily non-trivial.  
If $P$ is an isolated point in $LO(G)$ with associated ordering $<$ of $G$, then $P$ is certainly not an accumulation point of its conjugates in $LO(G)$.   It turns out that knowing $P$ is not an accumulation point of its conjugates $ gPg^{-1} \in LO(G)$ is enough to deduce that $B_<(G)$ (and hence $C_<(G)$) is non-trivial.

Observe that for any group $G$, if $1<h<g$ in the ordering corresponding to $P$, then left multiplication yields $1<h^{-1}g$, and then using the fact that $h$ is positive, we conclude that $1<h^{-1}gh$.  Translating this observation into a topological language, we have observed that if $P \in U_g$, then $hPh^{-1} \in U_g$ for any $h$ with $1<h<g$.  Supposing that
\[ \{P\} = \bigcap_{i=1}^m U_{g_i}, \]
is an isolated point, applying the above trick to the set of elements $\{g_1, \cdots , g_n \}$ allows us to conclude that for any $h$ with $1<h<g_i$ for all $i \in \{1, \cdots, n\}$, we must have 
\[ hPh^{-1} \in \bigcap_{i=1}^m U_{g_i}. \]
However, since $P$ is isolated, this means that $hPh^{-1} = P$, so that (in a sense soon to be made more precise) ``small elements in $G$ are bi-ordered,'' as they fix the positive cone $P$ under conjugation.

\begin{lemma}
\label{lem:subset}
Suppose that
\[ P \in \bigcap_{i=1}^m U_{g_i}, \]
where $\{ g_1, \cdots , g_m \}$  is some finite set of elements of
$G$, yet no conjugates of $P$ (different from $P$ itself) are in this open set.
 Then there exists $g_i \in \{ g_1, \cdots , g_m \}$ such
that the set

\[C_i = \{ g \in G : \mbox{ $g_i^{-k} \leq g \leq g_i^k$ for some
 $k$}\}\]
contains only elements of $G$ that fix the positive cone $P$ under
conjugation, that is, $g \in C_i \Rightarrow gPg^{-1}=P$.
\end{lemma}
\begin{proof} First, we show that there exists $g_i$ such that all elements in the set
\[ C_i^+ = \{ g \in G : \mbox{ $1<  g \leq g_i^k$ for some
 $k$}\}\]
fix $P$ under conjugation.

To this end, suppose not.  Then for each $g_i$ there exists $h_i$
with $1< h_i \leq g_i ^{k_i}$ for some $k_i$, and $h_iPh_i^{-1}
\neq P$. Choose $h = min\{h_1, \cdots , h_m\}$.  Then for each
$i$, we have
\[ h < g_i^{k_i} \Rightarrow  1< h^{-1}g_i^{k_i} \Rightarrow 1<
h^{-1}g_i^{k_i}h,\] and therefore $g_i^{k_i} \in hPh^{-1}$.  Now
since the element $g_i^{k_i}$ is positive in the order determined
by the positive cone $hPh^{-1}$, its $k_i$-th root $g_i$ is also
positive.  This shows that
\[ hPh^{-1} \in \bigcap_{i=1}^m U_{g_i}, \]
and by our choice of $h$, $hPh^{-1} \neq P$, a contradiction.
Therefore our claim holds for the set $C_i^+$.

To prove that all elements $g \in C_i$ fix the positive cone $P$, suppose that
$g \in G$ satisfies $g_i^{-k} \leq g <1 $ for some $k$.  Then $1
\leq g_i^k g < g_i^k$, so that either $g = g_i^{-k}$ or $g_i^k g
\in C_i^+$.
\begin{enumerate}
 \item In the case $g = g_i^{-k}$, then $g^{-1} \in C_i^+$ and so
fixes $P$, and so $g$ fixes $P$ under conjugation. \item If $g_i^k
g \in C_i^+$, then
\[ g_i^kgPg^{-1}g_i^{-k} = P,\]
so that we multiply by powers of $g_i$ from both sides and find
\[ gPg^{-1} = g_i^{-k}Pg_i^k =P.\]
Note that case (1) has been used to yield the final equality.
\end{enumerate}
Therefore we have found $g_i$ such that all elements in $C_i$ fix
$P$ as claimed.

\end{proof}

\begin{lemma}
\label{prop:convex} For any group $G$, if
\[ P \in \bigcap_{i=1}^m U_{g_i}, \]
and no conjugates of $P$ distinct from $P$ lie in this open set, then there exists
$g_i$ such that the set
\[C_i = \{ g \in G : \mbox{ $g_i^{-k} \leq g \leq g_i^k$ for some
 $k$}\}\]
 is a convex, bi-ordered subgroup of $G$.
\end{lemma}
\begin{proof}
Convexity of $C_i$ is clear from the definition.  By lemma
\ref{lem:subset}, $C_i$ is a subset of the bi-ordered group $stab(P)$, it follows that $C_i$ is bi-ordered
by the restriction ordering as well. Being bi-ordered, we can then
conclude that $C_i$ is a subgroup of $G$:  If $1<g \leq g_i^k$ for
some $k$, then $g_i^{-k} \leq g^{-1} < 1$, and similarly the implication $a<b \mbox{ and } c<d \Rightarrow ac < bd$ (this implication does not hold for left orders) shows closure under
multiplication.
\end{proof}

\begin{corollary}
\label{cor:BC}
Suppose that the left ordering $<$ of $G$ has positive cone $P$ which is not an accumulation point of its conjugates in $LO(G)$.  Then both $B_<(G)$ and $C_<(G)$ are non-trivial.
\end{corollary}

In particular, we have proven that if $<$ corresponds to an isolated point in $LO(G)$, then both $B_<(G)$ and $C_<(G)$ are non-trivial. 

We are now ready to complete the proof of Theorem \ref{thm:BCS}.
\begin{proof}[Proof of Theorem \ref{thm:BCS}.]
Let $P$ be the positive cone of a left ordering $<$ of a group $G$, and suppose that $P$ is an isolated point in $LO(G)$.  We know that $B_<(G)$ and $C_<(G)$ are non-trivial by Corollary \ref{cor:BC}, it remains to show that $B_<(G)$ is rank one abelian, and that $C_<(G)$ admits only finitely many left orderings.

Since we have assumed that $P$ is isolated in $LO(G)$, it follows from Lemma \ref{lem:congrp} that the restriction of $P$ to $B_<(G)$ must define a bi-ordering that is isolated in $LO(B_<(G))$.  However, by Theorem \ref{thm:bisolated}, this is only possible in the case when $B_<(G)$ is rank one abelian.   Similarly, It follows that the restriction of $P$ to $C_<(G)$ must define a Conradian ordering that is isolated in $LO(C_<(G))$, which by Theorem \ref{thm:consolated} is only possible in the case that $LO(C_<(G))$ is finite.
\end{proof}

\section{Dense and discrete orderings}

In recent work (\cite{CR07}, \cite{ALR08}), it has proven fruitful to consider discrete and dense group orderings separately, as they reflect different structures of the underlying group.  In considering the structure of $LO(G)$, dense orderings of a given group $G$ (with minor restrictions on the group $G$) are in some sense ``generic'' in $LO(G)$, in that dense orderings of $G$ constitute a dense $G_{\delta}$ set inside of a Cantor set within $LO(G)$.  Recall that a set $U$ in a topological space $X$ is a $G_{\delta}$ set if $U$ can be written as a countable intersection of open sets $\{ U_i \}_{i=1}^{\infty}$.

\begin{lemma}
\label{lem:gdelta}
Let $Z \subset LO(G)$ denote the set of dense left orderings of $G$.  If $G$ is countable, then $Z$ is a $G_{\delta}$ set.
\end{lemma}
\begin{proof}
Observe that if $\epsilon >1 $ is the least positive element in a left ordering $<$ of $G$ with positive cone $P$, then for all $g \in G$ (with $g \neq 1$ different from $\epsilon$) either $g< \epsilon^{-1}$ or $\epsilon <g$.  In other words, either $P \in U_{g^{-1}\epsilon^{-1}}$ or $P \in U_{\epsilon^{-1}g}$ for all $1 \neq g \in G$ different from $\epsilon$.  That is to say, let $V_{\epsilon}$ denote the set of all discrete left orderings of $G$ with least element $\epsilon$.  Then we have observed that
\[ V_{\epsilon}= \bigcap_{\epsilon \neq g \in G} (U_{g^{-1}\epsilon^{-1}} \cup  U_{\epsilon^{-1}g}) \cap U_{\epsilon}.
\] 
Note that $V_{\epsilon}$ is closed, as it is an intersection of closed sets, and consists of those positive cones that define an ordering of $G$ with $\epsilon$ as least positive element.  Therefore, the set of dense orderings is given by 
\[ Z = \bigcap_{1 \neq \epsilon \in G} (LO(G) \setminus V_{\epsilon}), \]
a countable intersection of open sets.
\end{proof}

The remaining difficulty is to show that any dense ordering is an accumulation point of other dense orderings. We first consider the case of abelian groups.

\subsection{Abelian groups}

From \cite{AK07}, we have the following fact:
\begin{proposition}
If $A$ is a torsion-free abelian group with $rank(A)>1$, then the space $LO(A)$
has no isolated points.
\end{proposition}

For a given torsion-free abelian group $A$, we can deduce much more about the structure of $LO(A)$ by examining the set of all dense orderings of $A$.

 \begin{proposition}
\label{prop:dense}
 Let $P$ be any positive cone in $LO(A)$, where $A$ is a torsion-free abelian group
 with $rank(A)>1$.  Then $P$ is an accumulation point of positive cones whose associated orderings are dense orderings.
 \end{proposition}

We begin by proving a special case.

\begin{lemma}
\label{lem:zk}
Let $P$ be any ordering in $LO( \mathbb{Z}^k )$, where $k>1$.  Then $P$ is an accumulation point of dense
 orderings.
\end{lemma}
\begin{proof}
We follow the ideas of Sikora in \cite{AS04}, making modifications where necessary.

For contradiction, let $k>1$ be the smallest $k$ for which the claim fails.  Suppose that 
\[ P \in \bigcap_{i=1}^n U_{g_i}, \]
with no dense orderings in this open set.  Note that we may assume that none of the $g_i$'s are integer multiples of one another.
Extend the ordering $<$ defined by $P$ to an ordering of $\mathbb{Q}^k$ by declaring $v_1<v_2$ for $v_1, v_2 \in \mathbb{Q}^k$ if $nv_1 <nv_2$ whenever $nv_1, nv_2 \in \mathbb{Z}^k$.  Let $H \subset \mathbb{Q}^k \otimes \mathbb{R} = \mathbb{R}^k$ be the subset of elements $x \in \mathbb{R}^k$ such that every Euclidean neighbourhood of $x$ contains both positive and negative elements.  Then $H$ is a hyperplane, and $H$ divides $\mathbb{R}^k$ into two components $H_-$ and $H_+$ having the property that $H_+$ contains only positive elements, and $H_-$ contains only negative elements.  Therefore the elements $g_i$ lie either in $H_+$ or $H$ itself.  
 
Suppose that two or more of the elements $\{g_1, \cdots , g_n \}$ lie inside $H$.  In this case, $H \cap \mathbb{Z}^k = \mathbb{Z}^m$ for some $m >1$ with $m<k$, and in this case the positive cone $P \cap \mathbb{Z}^m \subset \mathbb{Z}^m$ cannot be an accumulation point of dense orderings in $LO(\mathbb{Z}^m)$, for then we could change the positive cone $P$ using Lemma \ref{lem:congrp}.  This contradicts the minimality of $k$. 

The remaining possibilities are that exactly one (or none) of the elements  $\{g_1, \cdots , g_n \}$ lie inside $H$.  In this case, by slight perturbations of the hyperplane $H$, we can produce a new hyperplane $H'$ containing none of the lattice points $\mathbb{Z}^k \subset \mathbb{R}^k$, and with all points $g_i$ lying on one side of the hyperplane $H'$.  Specifically, if $\textbf{n}$ is the associated normal vector defining $H$, we may choose $H'$ having normal vector $\textbf{n}'$ arbitrarily close to $\textbf{n}$ (in the Euclidean distance), with the property that $\textbf{n}'$ has exactly one irrational entry.  This guarantees that no vector $\textbf{v} \in H'$ has all rational entries: If $\textbf{v}$ had all rational entries, the dot product $\textbf{n}' \cdot \textbf{v}$ would be a sum of $k-1$ rational numbers and one irrational number, and so cannot be zero.  Therefore, with normal vector $\textbf{n}'$ as above, $H' \cap \mathbb{Z}^k = \emptyset$.

This new hyperplane $H'$ defines a new ordering $P'$ on $\mathbb{Z}^k$ by declaring $P' =H_+' \cap \mathbb{Z}^k$, where $H_+'$ is the component of $\mathbb{R}^k \setminus H'$ containing all $g_i$.

To see that this ordering is dense, suppose that $\epsilon \in P'$ were a least element.  Then $\epsilon < \textbf{v}$ for all $\textbf{v} \in \mathbb{Z}^k$ iff $\textbf{v} - \epsilon \in H_+'$ for all $\textbf{v}$.  Consider the normal components $\epsilon_{\perp}$ and $\textbf{v}_{\perp}$ of $\epsilon$ and $\textbf{v}$.  As $H'$ contains no lattice points, $\epsilon_{\perp}>0$.  Then we can find $\textbf{v} \in P'$ with $\textbf{v}_{\perp} < \epsilon_{\perp}$, showing that $\textbf{v} - \epsilon \notin H_+'$ and $\epsilon$ is not the least positive element.
\end{proof}

\begin{proof}[Proof of Proposition \ref{prop:dense}]
To prove the statement for an arbitrary torsion-free abelian group $A$ with $rank(A)>1$, we let $g_1, \cdots ,g_m \in A$ be any finite family of elements in a given positive $P$.  We will show that there exist infinitely many positive cones with associated dense orderings on $A$ in which all $g_i$ are positive.  

Let $N$ be the subgroup of $A$ generated by the elements $g_1, \cdots, g_m$.  Then $N \cong \mathbb{Z}^k$ for $k \geq 1$.  Assume that $k>1$, for if it is the case that $N \cong \mathbb{Z}$, add an additional generator $g_{m+1}$ none of whose powers lie in $N$--we may do this since $rank(A)>1$.

By Lemma \ref{lem:zk}, $N$ admits infinitely many dense orderings in which all of $g_1, \cdots , g_m$ are positive, each constructed by perturbations of the hyperplane associated to the restriction order $P_N = N \cap P$.  Fix a positive cone $P_N'$ with a dense associated ordering of $N$, with $P_N' \neq P_N$.  We may extend $P_N'$ to a distinct ordering $Q$ on the isolator of $N$
\[I(N) = \{g \in A : \mbox{ $g^l \in N $ for some $l$} \} \] 
by declaring $g \in Q$ iff $g^l \in P_N'$ for some $l$.

Observe that the ordering of $I(N)$ with positive cone $Q$ is dense, for suppose not, say $Q$ had least element $\epsilon$.  Then $\epsilon \in P'$ is not possible since $P'$ is a dense ordering, so let $l>1$ be the least positive integer such that $\epsilon^l \in P'$.  By density of $P'$, we may then choose $g \in P' \subset Q$ with $1< g < \epsilon^l$.  Since the only positive elements less than $\epsilon^l$ are $\epsilon, \epsilon^2, \cdots, \epsilon^{l-1}$, we have that $g = \epsilon^i$ for $i<l$.  This contradicts our choice of $l$.

Now $I(N)$ is normal, and the quotient $A/I(N)$ is torsion-free abelian, so we may order the quotient. Using any ordering on the quotient, we can extend the dense ordering of $I(N)$ with positive cone $Q$ to give a dense ordering of $A$ with the required properties.
\end{proof}

Therefore, when $A$ is an abelian group with $rank(A)>1$, we know that the closure of the set of dense orderings in $LO(A)$ is the entire space $LO(A)$.  Thus, Proposition \ref{prop:dense} and Lemma \ref{lem:gdelta} together give us the following theorem.

\begin{theorem}
Suppose that $A$ is a countable abelian group.  Then $LO(A)$ is a Cantor set, and the set $Z$ of all dense left orderings of $A$ is a dense $G_{\delta}$ set within $LO(A)$.
\end{theorem}

Note that the case of discrete orderings must necessarily be different than this, for there exist abelian groups admitting no discrete orderings: divisible torsion free abelian groups are such an example.  Further, there exist abelian groups having a discrete set (in $LO(A)$) of discrete orderings.  As an example, consider $\mathbb{Z} \times \mathbb{Q}$.  This group has only four discrete orderings, namely the lexicographic orderings arising from the natural orderings (and their ``flipped'' versions) on both $\mathbb{Z}$ and $\mathbb{Q}$, with least positive elements $(1,0)$ and $(-1, 0)$ (Observe that  $(1,0)$ and $(-1, 0)$ are the only primitive elements in  $\mathbb{Z} \times \mathbb{Q}$, so any discrete ordering must have one of these elements as least positive element). 

\begin{question}
Let $A$ be a torsion-free abelian group with $rank(A)>1$.  What is the closure of the set of the discrete orderings in $LO(A)$?
\end{question}

\subsection{Non-abelian groups}

Our results concerning dense orderings generalize to the case of non-abelian groups.

\begin{proposition}
\label{prop:nonabdense}
Let $G$ be any group in which all rank one abelian subgroups are isomorphic to $\mathbb{Z}$.  If $P \in LO(G)$ corresponds to a dense left ordering $<$ of $G$, then $P$ is an accumulation point of positive cones whose associated left orderings are dense orderings.
\end{proposition}
\begin{proof}
Let $U=\bigcap_{i=1}^m U_{g_i}$ be an open set in $LO(G)$
containing $P$, the positive cone of a dense left ordering $<$ of $G$. 

If $U$ contains any conjugates of $P$ (different from $P$ itself),
then we are done, so suppose that no conjugate orderings lie in
$U$. Then by proposition \ref{prop:convex}, $G$ contains a convex,
bi-ordered subgroup $C$ of the form
\[C=C_i = \{ g \in G : g_i^{-k} \leq g \leq g_i^k \hspace{1em}
\mbox{for some
 $k$}\},\]
where $g_i \in \{g_1, \cdots , g_m\}$.  Denote by $C'$ the
intersection of all non-trivial convex subgroups of $C$.  There
are now two cases to consider.
\begin{enumerate}
\item $C' \neq \{1\}$.  In this case, since $C'$ is bi-ordered
and contains no convex subgroups, we can use a theorem of Conrad
which tells us the order must be Archimedian, and so $C'$ must be
abelian. From our assumption on $G$, if $rank(C')=1$, we have $C'
\cong \mathbb{Z}$, meaning our ordering is discrete.  Therefore
$rank(C')>1$.

Now the restriction ordering on $C'$ with positive cone $P \cap C'$ is a dense ordering, and we know from Theorem \ref{prop:dense} that every dense ordering in $LO(C')$ is an accumulation point of other dense orderings.  Therefore we may change the positive cone $P$ as in the proof of Lemma \ref{lem:congrp}, creating a new positive cone $P'$ containing all $g_i$, and corresponding to a dense ordering of $G$.

\item $C' = \{1 \}$.  In this case, $C$ must have infinitely
many convex subgroups whose intersection is trivial.  Therefore,
we may choose a convex subgroup $K$, that is non-trivial and
contains no $g_i$.  Define the positive cone of the ``flipped ordering'' of $K$ to be $(P^{-1} \cap K) = P_K^{-1}$.  Then we define a new positive cone $P' \subset G$, with $P' \in U$, by setting $P' = P_K^{-1} \cup (P \cap G \setminus K)$.  Again, the new ordering $<'$ of $K$ with positive cone $P'$ is dense, and so the ordering we have defined on $G$ is dense.
\end{enumerate}
\end{proof}

In the case of an abelian group $A$, the closure of the set of dense orderings was the entire space $LO(A)$, which is known to be homeomorphic to the Cantor set when $A$ is countable.   In the non-abelian case, Theorem \ref{th:dense} gives us a similar result.

\begin{proof} [Proof of Theorem \ref{th:dense}]
Let $G$ be any countable group with all rank one abelian subgroups isomorphic to $\mathbb{Z}$.  Then since $G$ is countable, $LO(G)$ is metrizable, as is the space $\bar{Z} \subset LO(G)$.  Proposition \ref{prop:nonabdense} shows that the set $\bar{Z}$ contains no isolated points, and since it is closed, it is compact.  Therefore $\bar{Z}$ is a compact, metrizable, totally disconnected perfect space, and so is homeomorphic to the Cantor set \cite{HY61}.  By Lemma \ref{lem:gdelta}, the set $Z$ is also a dense $G_{\delta}$ set within $\bar{Z}$.
\end{proof}
 
With the restriction that all rank one abelian subgroups of $G$ be isomorphic to $\mathbb{Z}$, it also follows readily that any isolated point in $LO(G)$ must correspond to a discrete left ordering of $G$.  This can be seen by appealing to either Theorem \ref{th:dense} (which is stronger than what we need), or by appealing to Theorem \ref{thm:BCS}, and remarking that the smallest convex subgroup in the Conradian soul of an isolated left ordering must be a rank one abelian group.  

We turn our attention next to discrete orderings, and observe conditions under which a discrete ordering of $G$ is not an isolated point in $LO(G)$.  We no longer need the restriction that all rank one abelian subgroups be isomorphic to $\mathbb{Z}$.

\begin{lemma}
\label{lem:condisc}  Suppose that $P \subset G$ is the positive cone of a discrete left ordering $<$ with least element $\epsilon$.  Then if
$g \epsilon g^{-1}> 1$ for all $g \in G$ and
\[ P \in \bigcap_{i=1}^m U_{g_i} \]
contains no conjugates of $P$, there exists $g_i$ which is not a
power of $\epsilon$  such that
\[C_i = \{ g \in G : \mbox{$g_i^{-k} \leq g \leq g_i^k$ for some
 $k$}\}\]
 is a convex, bi-ordered subgroup which properly contains the convex subgroup $\langle \epsilon \rangle$.
\end{lemma}
\begin{proof}
Suppose that $U = \bigcap_{i=1}^m U_{g_i}$ contains $P$, but no
conjugates of $P$.  If no $g_i$ is equal to a power of $\epsilon$,
then we are done, as we may apply proposition
\ref{prop:convex}.

On the other hand, suppose that some $g_i$ is a power of $\epsilon$, say $g_1 = \epsilon^l$.  Then the
condition $g \epsilon g^{-1} > 1 $ for all $g \in G$ guarantees that the
open set $U_{\epsilon}$ contains every conjugate of $P$.
Therefore, if  
\[ \left( \bigcap_{i=2}^m U_{g_i} \right)  \cap U_{\epsilon^l} \]
contains no conjugates of $P$, neither does the open set
$\bigcap_{i=2}^m U_{g_i}$.  Continuing to eliminate powers of
$\epsilon$ in this way, we can eventually find an open set
$\bigcap_{i=r}^m U_{g_i}$ containing no conjugates of $P$, and
with no $g_i$ equal to a power of $\epsilon$.  From here we may
apply Proposition \ref{prop:convex}.
\end{proof}

\begin{theorem}
\label{th:discacc}
Let $G$ be a group, and $P$ the positive cone of a discrete left ordering $<$ with least positive element $\epsilon$.  If $g \epsilon g^{-1} \in P$ for all $g \in G$, then $P$ is not isolated in $LO(G)$.
\end{theorem}
\begin{proof}
We proceed very similarly to the proof of Theorem \ref{th:dense}.  Let $U=\bigcap_{i=1}^m U_{g_i}$ be an open set in $LO(G)$ containing $P$.  If $U$ contains any conjugates of $P$, then we are done, so by Lemma \ref{lem:condisc}, we may suppose that there exists convex subgroup $C$ properly containing $\left< \epsilon \right>$, which is bi-ordered by the restriction of $P$.  

Note that the convex subgroup $C$ is not rank one abelian:  Suppose that $rank(C)=1$.  As the containment $\left< \epsilon \right> \subset C$ is proper, we can choose $c \in C$ with $C \neq 1$,  that is not a power of $\epsilon$.  If we then assume that $C$ is rank one abelian, we arrive at $\epsilon^k = c^l$ for some integers $k, l$, contradicting the fact that $\epsilon$ is the least positive element.

Thus, by Theorem \ref{thm:bisolated}, we know that the restriction of $P$ to the subgroup $C$ is not isolated in $LO(C)$, and it follows from Lemma \ref{lem:congrp} that $P$ is not isolated in $LO(G)$.
\end{proof}

\section{The braid groups}

As a sample application of these results, we turn our focus to the braid groups.  It is known that the space of left orders $LO(B_n)$ is \textit{not} homeomorphic to the Cantor set for $n \geq 2$.  We begin by defining the Dehornoy left ordering of the braid groups (also known as the `standard' ordering), whose positive cone we shall denote $P_D$ \cite{DDRW08}, \cite{PD94}. Recall that for each integer $n \ge 2$, the Artin braid group $B_n$ is the group generated by 
$\s_1 , \s_2 , \dots , \s_{n-1}$, subject to the relations 
$$\s_i\s_j = \s_j\s_i {\rm \: if \:} |i-j| >1,\quad \s_i\s_j\s_i=\s_j\s_i\s_j {\rm \: if \: } |i-j| =1.$$  

\begin{definition}Let $w$ be a word in the generators $\s_i, \cdots , \s_{n-1}$.  Then $w$ is said to be: $i$-positive if the generator $\s_i$ occurs in $W$ with only positive exponents, $i$-negative if $\s_i$ occurs with only negative exponents, and $i$-neutral if $\s_i$ does not occur in $w$.
\end{definition}

It is shown in \cite{PD94} that for every integer $i$ with $1 \leq i < n$, every braid $\beta \in B_n$ is either $i$-positive, $i$-negative, or $i$-neutral. We may then define the positive cone of the Dehornoy ordering as

\begin{definition} 
The positive cone $P_D \subset B_n$ of the Dehornoy ordering is the set
\[P_D = \{ \beta \in B_n : \mbox{ $\beta$ is $i$-positive for some $i \leq n-1$}\}.\]
\end{definition}

There is also a second positive cone of interest, discovered by the authors of \cite{DD04}, which we shall denote by $P_{DD}$.  Denote by $P_i \subset B_n$ the set of all $i$-positive braids.  Note that the set of all $i$-negative braids is simply $P_i^{-1}$.

\begin{definition} The positive cone $P_{DD} \subset B_n$ is the set
\[P_{DD} = P_{n-1} \cup P_{n-2}^{-1} \cup \cdots \cup P_{1}^{(-1)^{n}}.\]
\end{definition}

That either of these notions defines a positive cone in $B_n$ is difficult to show, as it is not clear that the notion of a braid being $i$-positive is well defined.  This was the main idea introduced to braid theorists in Dehornoy's seminal paper \cite{PD94}. 

The positive cone $P_{DD}$ was originally defined in light of the following property:
\begin{proposition} [Dubrovina, Dubrovin \cite{DD04}] The positive cone $P_{DD}$ is generated as a semigroup by the braids
\[ y_1 = \s_1 \cdots \s_{n-1}, y_2 = (\s_2 \cdots \s_{n-1})^{-1}, y_3 = \s_3  \cdots \s_{n-1}, \cdots,  y_{n-1}=\s_{n-1}^{(-1)^n}. \]
\end{proposition}

Note that for two positive cones $P$ and $Q$, if $P\subset Q$ then necessarily $P=Q$.  Therefore
\begin{corollary}
The order $P_{DD}$ is an isolated point in $LO(B_n)$, in particular, 
\[ \{P_{DD}\} = \bigcap_{i=1}^{n-1} U_{y_i} .\]
\end{corollary}

Knowing that $LO(B_n)$ has isolated points for $n \geq 2$, it makes sense to ask the question:  Is the standard ordering $P_D$ an isolated point in $LO(B_n)$?  This question is answered in \cite{DDRW08}, using a very explicit calculation.  That $P_D$ is not isolated, however, was originally proven in \cite{NF07}, though the techniques are different than those used here, which illustrate our machinery.

First, we begin with a proposition which establishes a very important property of the ordering $P_D$.  Recall the Garside monoid $B_n^+ \subset B_n$ is the monoid generated by the elements $\s_1, \cdots \s_n$.
\begin{proposition}
Let $\beta \in B_n$ and $\alpha \in B_n^+$ be given.  Then $\beta \alpha \beta^{-1} \in P_D$.
\end{proposition}
This property of the Dehornoy ordering is referred to as the subword property, or property \textbf{S}. 

Next, we must know that the Dehornoy ordering is discrete \cite{CR07}.  

\begin{proposition}
The Dehornoy ordering of $B_n$ is discrete, with smallest positive element $\s_{n-1}$.  
\end{proposition}

These two propositions together show us that $P_D$ satisfies the hypotheses of Theorem \ref{th:discacc}.  If we can additionally show that $P_D$ has no biorderable convex subgroups properly containing $\left< \s_{n-1} \right>$, then we can conclude that $P_D$ is an accumulation point of its orbit under the $B_n$-action on $LO(B_n)$.

Recall the natural inclusions 
$B_m \subset B_n$ whenever $m \le n$ which takes $\s_i \in B_m$ to $\s_i \in B_n$.
A useful operation is the shift homomorphism $sh \colon B_m \to B_n, \; m < n$ defined by $sh(\s_i) = \s_{i+1}$.  This is clearly injective and order-preserving.  The shift may be iterated, and we note that 
$sh^r(B_{n-r})$ is just the subgroup $\langle \s_{r+1}, \dots \s_{n-1} \rangle$ of $B_n$, or in other words, the subgroup of all elements which are $i$-neutral for all $i \le r$.  

\begin{lemma}
The subgroups $sh^r(B_{n-r})$, $r>0$, are the only convex subgroups under the ordering $P_D$. 
\end{lemma}
\begin{proof}
Set $H_r = sh^r(B_{n-r})$, and let $C$ be a convex subgroup in the Dehornoy ordering.  Choose $i$ to be the smallest integer such that $C$ contains an $i$-positive braid.  Then clearly $C \subset H_{i+1}$, our aim is to show the opposite inclusion, which establishes the claim. 

Let $\beta \in C$ be an $i$-positive braid.  The braid $\s_j^{-1} \beta$ is $i$-positive for $j>i$, so that $1 < \s_j < \beta \Rightarrow \s_j \in C$, and so $H_{i+1} \subset C$.  Considering the generator $\s_i$, we write $\beta = w_1 \s_i w_2$, where $w_1$ is an empty or $i$-neutral word, and $w_2$ is an empty, $i$-neutral, or $i$-positive word.  We will show $\s_i \in C$.

First, we note that the the braid represented by the word $\s_i w_2$ lies in $C$, as $w_1$ contains only $\s_{i+1}, \cdots , \s_{n-1}$, all of which are in $C$.  If $w_2$ is empty, the claim is proven, if $w_2$ is $i$-neutral, then we may right multiply by appropriate $\s_j$ for $j>i$ to arrive at $\s_i \in C $, and again the claim is proven.  Lastly, if $w_2$ is $i$-positive, then we get:
\[  1< w_2 \Rightarrow 1<\s_i < \s_i w_2 \in C,\]
and the claim follows from convexity of $C$.
\end{proof}

Since all convex subgroups are isomorphic to a shifted copy of the braid groups, we conclude that
\begin{corollary}
\label{cor:deho}
No subgroup that is convex under the ordering $P_D$ is bi-orderable, except for the subgroup $\left< \s_{n-1} \right>$.
\end{corollary}
 
\begin{theorem}
For every $n>2$, the positive cone $P_D$ in $B_n$ is an accumulation point of its conjugates in $LO(B_n)$.
\end{theorem}
\begin{proof}
Apply Corollary \ref{cor:deho} and Lemma \ref{lem:condisc}.
\end{proof}

\bibliographystyle{plain}
\bibliography{candidacy}

\end{document}